\documentclass[12pt]{article}
\usepackage{amsmath,amssymb,graphicx,bm}
\usepackage{graphicx}
\usepackage{epsfig}

\evensidemargin=\oddsidemargin 
\newdimen\dummy
\dummy=\oddsidemargin 
\newtheorem{theorem}{Theorem}[section]
\newtheorem{lemma}[theorem]{Lemma}

\newcommand{\be}{\begin{eqnarray}}
\newcommand{\ee}{\end{eqnarray}}
\newcommand{\ben}{\begin{eqnarray*}}
\newcommand{\een}{\end{eqnarray*}}
\newcommand{\bq}{\begin{equation}}
\newcommand{\eq}{\end{equation}}
\newcommand{\bqn}{\begin{equation*}}
\newcommand{\eqn}{\end{equation*}}

\begin {document}

\title{\large \bf The Heterogeneous Multiscale Finite Volume
Method for Convection-Diffusion-Reaction Problem}

\author{ Tao Yu\footnote{Department of Mathematics and Physics,
Jinggangshan University, Ji'an 343009, China.
(yutao@jgsu.edu.cn)} \and  Haitao Cao\footnote{Department of Mathematics and Physics, Hohai University, Changzhou Campus, Changzhou 213022, China.
(20104007004@suda.edu.cn)} }
\date{}

\maketitle

\begin{abstract}
In this paper, we employ an finite volume method (FVM) based on the heterogenous
multiscale method (HMM), for the multiscale convection-diffusion-reaction
problem. The optimal order convergence rate in $H^1$-norm is given for periodic medias.
\end{abstract}


\section{Introduction}

This paper consider the multiscale method for the following convection-diffusion-reaction problem
\be\label{origin
problem}\left\{
\begin{array}{ll}\displaystyle
-\nabla(a^{\varepsilon}(x)\nabla
u^{\varepsilon}(x))+b^{\varepsilon}(x)\nabla u^{\varepsilon}(x)+c^{\varepsilon}(x)u^{\varepsilon}(x))
=f(x),\;\;&x\in\Omega,\\\displaystyle
u^{\varepsilon}(x)=0,&x\in\partial\Omega,\end{array}\right. \ee
where $\Omega\subset\mathbb{R}^2$ (or $\mathbb{R}^3$) is a bounded convex polygonal
domain with a Lipschitz boundary $\partial \Omega$. $\varepsilon\ll 1$ is a positive parameter which signifies the
multiscale nature of the problem. This problem is related to the studying of groundwater and solute transport in porous media(see \cite{Bear}).

Optimal order convergence rate of classical finite element method based on piecewise linear polynomials relies on the $H^2$-norm of $u^{\varepsilon}$. As the coefficient varies on a scale of $\varepsilon$, the solution $u^{\varepsilon}$ may also  oscillate at the same scale. A direct numerical solution of this multiscale problem is very difficult unless the mesh size is smaller enough. However, this is not computationally feasible in many applications. On the other hand, from an engineering perspective, the macroscopic features of the solution are often of the main interest and importance. Through the homogenization theory \cite{Bemsoussan,Jikov}, there is a homogenized equation which can capture the macroscopic properties. That is to say there exists homogenized coefficients $a^0$, $b^0$, $c^0$, such that
\be\label{homogenized
problem}\left\{
\begin{array}{ll}\displaystyle
-\nabla(a^{0}(x)\nabla
u^{0}(x))+b^{0}(x)\nabla u^{0}(x)+c^{0}(x)u^{0}(x))
=f(x),\;\;&x\in\Omega,\\\displaystyle
u^0(x)=0,&x\in\partial\Omega.\end{array}\right. \ee
It is obvious
that the classical methods are effective to solve this equation (\ref{homogenized problem}). Unfortunately, in general there are no explicit formulas for the homogenized coefficients, except the restrictive assumptions on the media. Therefore, it is desirable to develop numerical methods that can capture the effect of small scale on the large scale. To overcome the difficulty, many methods are designed to solve the problems on grids which are coarser than the scale of oscillation, see, for example, \cite{Babuska}, \cite{E1}, \cite{Franca}, \cite{Hou1}, \cite{E1}, \cite{Hughes} and references therein.

The heterogeneous multiscale method (HMM) was introduced in \cite{E1}. This method is a general efficient methodology for the multiscale problems. It consists of two components: selecting a macro-scopic solver on a coarse mesh and estimating the missing macroscale data by solving the local fine scale problems. The careful choosing of the macro solver and local fine problems is the key issue for this method. A different choice of macro-scopic solver leads to a different heterogeneous multiscale method. In this paper we choose the finite volume methods (FVM) introduced in \cite{Wu} as the macro-scopic solver. For convenience, the heterogeneous multiscale method taking the finite volume method as the macroscopic solver is called HMM-FVM for abbreviation.

In the remained part of this paper, one HMM-FVM for the multiscale convection-diffusion-reaction problem is present in Section 2, then the approximate solution and its error estimates for periodic media are shown in Section 3.
\section{HMM-FVM for convection-diffusion-reaction problem}
\setcounter{equation}{0}

In this section, we first introduce the finite volume method for convection-diffusion-reaction problem in \cite{Wu}. Then by taking this method as a macro solver of HMM, we derive the HMM-FVM.

Let $T_H$ be a quasi-uniform triangulation of the polygonal domain $\Omega$. The barycenter dual decomposition $T_H^*$ is constructed by connecting the barycenter to the midpoints of edges of every triangle element by straight lines. Given a triangle element $K$, let $\tau$ and $Q$ be an edge and barycenter, respectively, of $K$, $R$ be the midpoint of $\tau$. Next, let $P$ be a modal point and $K_P^*$ be the dual element with respect to $P$ (see Figure 1). Denote by $H$ the maximum length of the edges, $N_H$ the set of all nodal points. And denote by $\dot{K}$ the vertices of $K$.

\begin{figure}[htbp]
\centering
 \scalebox{0.5}[0.5]{\includegraphics{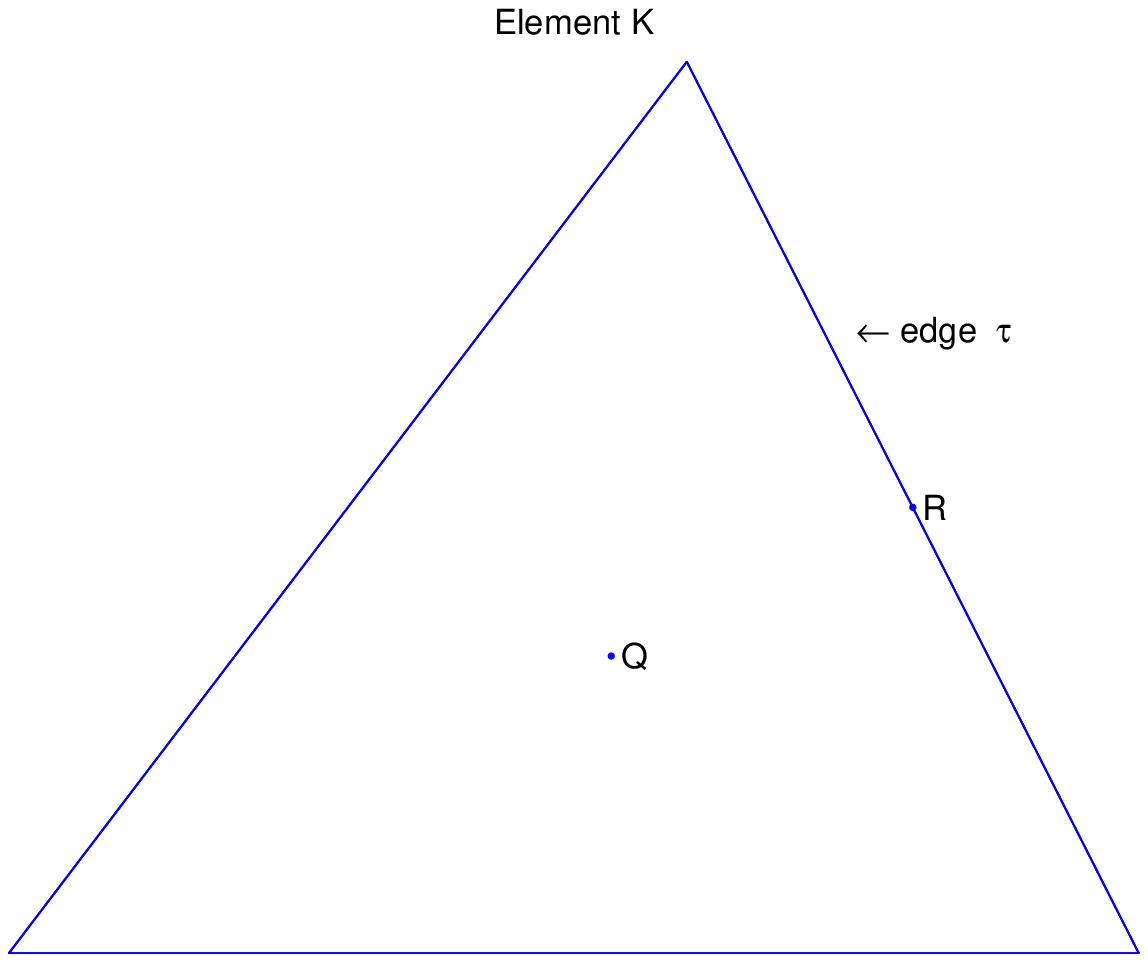}}
 \scalebox{0.5}[0.5]{\includegraphics{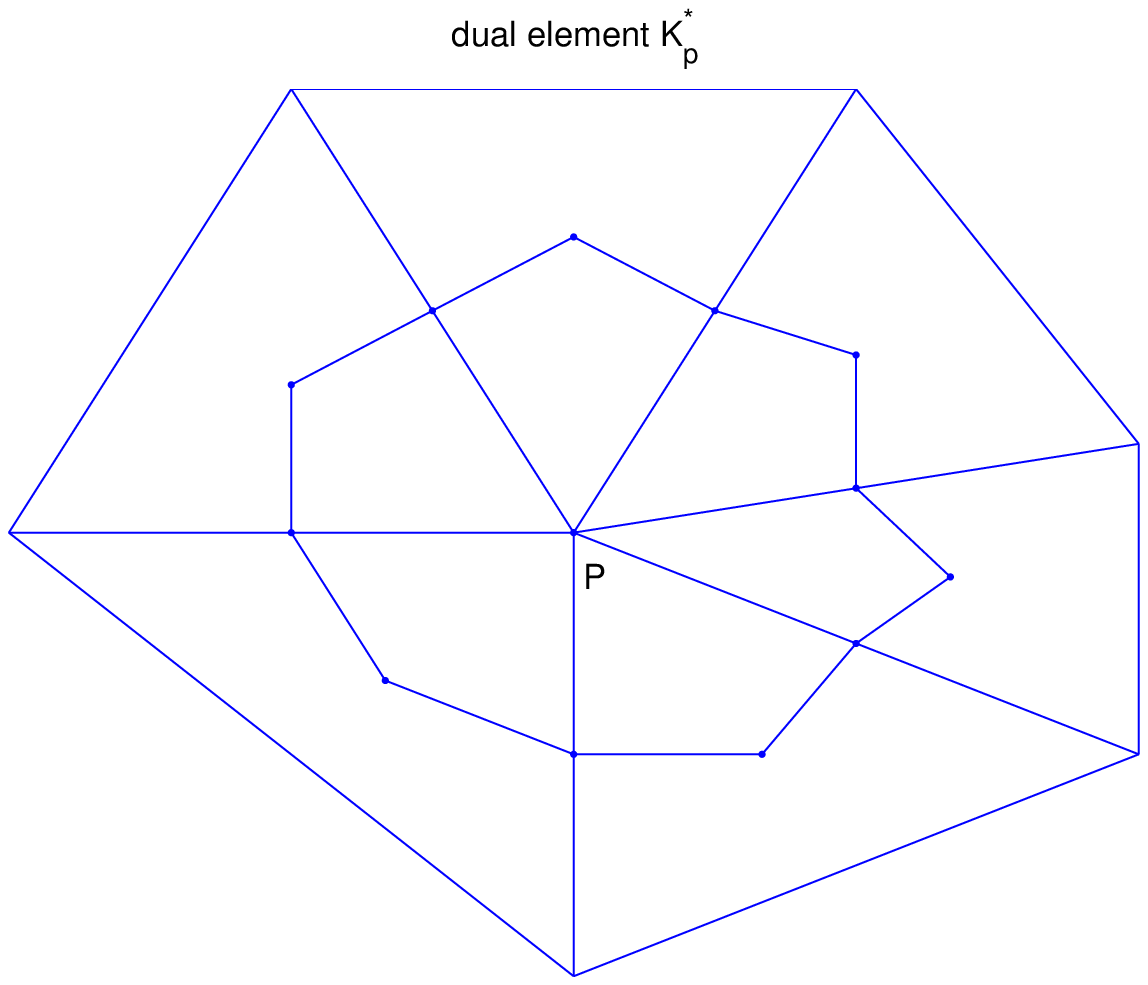}}
 \caption{The left figure is an element K, and the right figure represents the dual element with
respect to a node P}\label{fig1}
\end{figure}

Let $U=H_0^1(\Omega)=\{u\in H^{1}(\Omega)| u=0\;on\;\partial\Omega\}$ and $U_H \subset U$ be the piecewise linear finite element space on $T_H$. Denote by $ \Pi_H^*$ the interpolation operator from $U_H$ to the piecewise constant space on $T_H^*$:
$$(\Pi_H^*v_H)(x)=v_H(P), x\in K_P^*, \forall P\in N_H, v_H\in U_H.$$
Then the finite volume method(FVM)\cite{Wu} reads as: Find $u_H\in U_H$ such that
\be A(u_H, \Pi_H^* v_H)=(f, \Pi_H^* v_H), \forall v_H\in U_H,\ee
where
\be\label{A}
A(u_H, \Pi_H^* v_H)=\sum_{P\in N_H}\left[-\int_{\partial K_P^*}\mathbf{n}\cdot(a^H\nabla u_H)\Pi_H^* v_H ds
+\int_{K_P^*}(b^H\nabla u_H+c^H u_H)\Pi_H^* v_H dx\right]\nonumber\\
=\sum_{K\in T_H}\left[-\sum_{P\in\dot{K}}\int_{K\cap\partial K_P^*}\mathbf{n}\cdot(a^H\nabla u_H)\Pi_H^* v_H ds
+\int_{K}(b^H\nabla u_H+c^H u_H)\Pi_H^* v_H dx\right],
\ee
and
\be\label{f1}
(f, \Pi_H^* v_H)=\sum_{K\in T_H}\int_{K}f \Pi_H^* v_H dx.
\ee

Consider the numerical integration for the barycenter quadrature rule, (\ref{A}) and (\ref{f1}) can be written as:
\be
A_{FVM}(u_H, \Pi_H^* v_H)=\sum_{K\in T_H}[-\sum_{P\in\dot{K}}\left|K\cap\partial K_P^*\right|\mathbf{n}\cdot(a^H(Q)\nabla u_H)\Pi_H^* v_H\nonumber \\
+\left|K\right|\left(b^H(Q)\nabla u_H+c^H(Q) u_H(Q)\right)\Pi_H^* v_H ],
\ee
and
\be
(f, \Pi_H^* v_H)_H=\sum_{K\in T_H}\left|K\right|f(Q) \Pi_H^* v_H.
\ee
Thus, the barycenter quadrature approximation of the FVM reads: Find $u_H \in U_H$ such that
\be\label{HMM_FVM}
A_{FVM}(u_H, \Pi_H^* v_H)=(f, \Pi_H^* v_H)_H,\;\;\;\forall v_H \in U_H.
\ee

This paper considers the scale separation in the coefficients. Assume the multiscale coefficients $a^{\varepsilon}(x)$, $b^{\varepsilon}(x)$ and $c^{\varepsilon}(x)$ have the forms $a(x, \frac{x}{\varepsilon})$, $b(x, \frac{x}{\varepsilon})$ and $c(x, \frac{x}{\varepsilon})$. Moreover, assume that $a(x, y)$, $b(x, y)$ and $c(x, y)$ are smooth in $x$ and periodic in $y$ with respect to the unit cube $Y$.

In the absence of explicit knowledge of $a^H$, $b^H$ and $c^H$, $a^H(Q)\nabla u_H$, $b^H(Q)\nabla u_H$ and $c^H(Q)$ can be approximated by
\be\label{effective a}
a^H(Q)\nabla u_H\simeq \frac{1}{|K_{\delta}(Q)|}\int_{K_{\delta}(Q)}a(Q,
\frac{x}{\varepsilon})\nabla R(u_H)dx,
\ee
\be\label{effective b}
b^H(Q)\nabla u_H\simeq \frac{1}{|K_{\delta}(Q)|}\int_{K_{\delta}(Q)}b(Q,
\frac{x}{\varepsilon})\nabla R(u_H)dx,
\ee
and
\be\label{effective c}
c^H(Q) \simeq \frac{1}{|K_{\delta}(Q)|}\int_{K_{\delta}(Q)}c(Q,
\frac{x}{\varepsilon})dx,
\ee
 where $R(u_H)$ is the solution of the following microcell problem
\bq\label{micro}
\left\{\begin{array}{cc}\displaystyle -\nabla\cdot(a(Q,
\frac{x}{\varepsilon})\nabla R(u_H))=0,\;\;\;
&\mbox{in}\;K_{\delta}(Q),\\\displaystyle
R(u_H)=u_H,&\mbox{on}\;
\partial
K_{\delta}(Q),\end{array}\right.\eq
and $K_{\delta}(Q)$ is a cube of size $\delta$ centered at $Q$.

Thus, the HMM-FVM is defined by the following variational problem: $\forall v_H\in U_H$, find $u_H\in U_H$  that satisfies (\ref{HMM_FVM}), (\ref{effective a}), (\ref{effective b}) and (\ref{effective c}).

Our main result is the following.
\begin{theorem}\label{error estimate} Denote by $u^0$ and $u_H$
the solution of (\ref{homogenized problem}) and (\ref{HMM_FVM}),
respectively. Then for sufficient small $H$ and
$\varepsilon/\delta$, we have \ben||u^0-u_H||_{1, \Omega}\leq C\left(H+\frac{\varepsilon}{\delta}\right),\een where $C$ is a positive constant independent of $\varepsilon$, $\delta$ and $H$.
\end{theorem}

\section{A priori error estimate}
\setcounter{equation}{0}

In this section, we consider the HMM-FVM as a perturbation of the linear finite element method(FEM), and then obtain the optimal $H^1$ error estimate.

The finite element method for the homogenization problem (\ref{homogenized problem}) reads as:
find $u_H\in U_H$ such that
\be\label{FEM}
B(u_H, v_H)=(f, v_H),~~~\forall v_H\in U_H,
\ee
where
\be\label{B}
B(u_H, v_H)=\sum_{K\in T_H}\int_{K}\left(a^0\nabla u_H\cdot\nabla v_H+b^0\nabla u_Hv_H+c^0 u_H v_H\right)dx
\ee
and
\be\label{f}
(f, v_H)=\sum_{K\in T_H}\int_{K}f v_H dx.
\ee

Using the barycenter quadrature into (\ref{B}), the  bilinear form can be refined by
\be\label{FEM int}
\tilde{B}(u_H, v_H)=\sum_{K\in T_H}|K|\left(a^0(Q)\nabla u_H\cdot\nabla v_H+b^0(Q)\nabla u_Hv_H(Q)+c^0(Q) u_H(Q) v_H(Q)\right).
\ee

The following Lemma characterizes the difference between the bilinear form of the HMM-FVM and that of the FEM, which plays the key role in the subsequent analysis.
\begin{lemma}\label{B_A}
$\forall u_H, v_H\in U_H$, we have
\ben
|B(u_H, v_H)-A_{FVM}(u_H, v_H)|\leq C \left( H+\frac{\varepsilon}{\delta}\right)||u_H||_{1,~\Omega}||v_H||_{1,~\Omega},
\een
where $C$ is a positive constant independent of $\varepsilon$, $\delta$ and $H$.
\end{lemma}
Proof. The bilinear form can be eliminated by three terms
\ben
|B(u_H, v_H)-A_{FVM}(u_H, v_H)|&\leq&|B(u_H, v_H)-\tilde{B}(u_H, v_H)|\\
&+&|\tilde{B}(u_H, v_H)-\tilde{A}(u_H, v_H)|\\
&+&|\tilde{A}(u_H, v_H)-A_{FVM}(u_H, v_H)|\\
&:=&\varepsilon_1(u_H, v_H)+\varepsilon_2(u_H, v_H)+\varepsilon_3(u_H, v_H),
\een
where
\be
\tilde{A}(u_H, v_H)=\sum_{K\in T_H}|K|\left(a^H(Q)\nabla u_H\cdot\nabla v_H+b^H(Q)\nabla u_Hv_H(Q)+c^H(Q) u_H(Q) v_H(Q)\right).
\ee

Denote by
\be
E_K(f)=\int_K f(x)dx-|K|f(Q),
\ee
then from the standard estimate\cite{Clarlet}
\be
|E_K(a(x)p(x)q(x))|\leq CH||a||_{1,~\infty}||p||_{0,~K}||q||_{0,~K},
\ee
the numerical quadrature error $\varepsilon_1(u_H, v_H)$ can be easily estimated as
\be
\varepsilon_1(u_H, v_H)&\leq& \sum_{K\in T_H}\left(|E_K(a^0\nabla u_H\cdot \nabla v_H)|
+|E_K(b^0u_H \nabla v_H)|+|E_K(c^0u_H v_H)|\right)\\
&\leq& CH||u_H||_{1,~\Omega}||v_H||_{1,~\Omega}.
\ee

In order to estimate the modeling error $\varepsilon_2(u_H, v_H)$, the following Lemmas are useful.
\begin{lemma}\label{model error for a and b}(\cite{Chen, E2, Hou})
There exists a constant $C$ independent of H, $\varepsilon$ and
$\delta$ such that for $a^H$, $b^H$,
\ben \left|a^H_{ij}-a^0_{ij}\right|\leq
C\frac{\varepsilon}{\delta},\;\;i, j=1,2,\een
and
\ben\left|b^H_{i}-b^0_{i}\right|\leq
C\frac{\varepsilon}{\delta},\;\;i=1,2.\een
\end{lemma}

\begin{lemma}\label{Lemma c}(\cite{Chen, Chen2})
Given domain $K_{\delta}$ with $diam(K)=\delta$, let $\varphi(s, y)$ defined in $Y$ be a $Y-$periodic function in $y$, where $Y$ is a unit cube and $s\in R$ is fixed. Then
\ben \left|\frac{1}{|Y|}\int_Y \varphi(s, y)dy-<\varphi(s, x/\varepsilon)>_K\right|\leq
C\frac{\varepsilon}{\delta},\een
where $C$ is independent of $\varepsilon$, $\delta$ and $s$.
\end{lemma}

Directly from the definition of $c^H$ in (\ref{effective c}) and Lemma \ref{Lemma c}, the error between $c^H$ and $c^0$ can be easily estimate as
\be\label{ch_c0} \left|c^H-c^0\right|\leq
C\frac{\varepsilon}{\delta},\ee
where $C$ is a positive constant independent of $\varepsilon$, $\delta$ and $H$.

Denote by $$e(HMM)=\max_{x\in \Omega}\left\{\left|a^H_{ij}-a^0_{ij}\right|, \left|b^H_{i}-b^0_{i}\right|, \left|c^H-c^0\right|\right\},$$
then
\be\label{eHMM}
e(HMM)\leq C\frac{\varepsilon}{\delta}
\ee
from Lemma \ref{model error for a and b} and (\ref{ch_c0}).\\
So the modeling error $\varepsilon_2(u_H, v_H)$ can be estimate as
\be
\varepsilon_2(u_H, v_H)\leq C e(HMM)||u_H||_{1,~\Omega}||v_H||_{1,~\Omega}\leq C\frac{\varepsilon}{\delta}||u_H||_{1,~\Omega}||v_H||_{1,~\Omega}.
\ee

It remains to estimate the term $\varepsilon_3(u_H, v_H)$. Deduce by Green's formula, the first term of $\tilde{A}(u_H, v_H)$ reads
\ben
|K|a^{H}(Q)\nabla u_{H}\nabla v_{H}&=&\int_{K}a^{H}(Q)\nabla u_{H}\nabla v_{H}dx\\
&=&-\int_{K}\nabla(a^{H}(Q)\nabla u_{H})v_H dx+\int_{\partial K}n\cdot(a^{H}(Q)\nabla u_{H})v_Hds\\
&=&\int_{\partial K}n\cdot(a^{H}(Q)\nabla u_{H})v_Hds
\een
The first term of $A_{FVM}(u_H, v_H)$ can derived similarly
\ben
&-&\sum_{P\in \dot{K}}|K\cap\partial K_P^*|n\cdot(a^{H}(Q)\nabla u_{H})\Pi_H^*v_H\\
&=&-\sum_{P\in \dot{K}}\int_{K\cap\partial K_P^*}n\cdot(a^{H}(Q)\nabla u_{H})\Pi_H^*v_Hds\\
&=&-\sum_{P\in \dot{K}}\int_{\partial(K\cap\partial K_P^*)}n\cdot(a^{H}(Q)\nabla u_{H})\Pi_H^*v_Hds+\int_{\partial K}n\cdot(a^{H}(Q)\nabla u_{H})\Pi_H^*v_Hds\\
&=&-\int_{K}\nabla(a^{H}(Q)\nabla u_{H})\Pi_H^*v_Hdx+\int_{\partial K}n\cdot(a^{H}(Q)\nabla u_{H})\Pi_H^*v_Hds\\
&=&\int_{\partial K}n\cdot(a^{H}(Q)\nabla u_{H})\Pi_H^*v_Hds
\een
Thus by (\ref{eHMM}) and the property of the operator $\Pi_H^*$ (see Lemma 2.1 in \cite{Wu})
\be
||v_H-\Pi_H^*v_H||_{0,K}\leq CH|v_H|_{1,K},\;\;\forall v_H \in U_H,
\ee
the last term $\varepsilon_3(u_H, v_H)$ can be estimated as
\ben
\varepsilon_3(u_H, v_H)\!\!\!&=&\!\!\!\sum_{K\in T_H}[\sum_{e\in \partial K}\int_{e} n\cdot(a^{H}(Q)\nabla u_{H})(v_H-\Pi_H^*v_H) ds \\
\!\!\!&&\!\!\!+|K|\left(b^{H}(Q)\nabla u_{H}+c^{H}(Q) u_{H}\right)(v_H-\Pi_H^*v_H)]\\
\!\!\!&=&\!\!\!\sum_{K\in T_H}|K|\left(b^{H}(Q)\nabla u_{H}+c^{H}(Q) u_{H}\right)(v_H-\Pi_H^*v_H)\\
\!\!\!&\leq &\!\!\!\sum_{K\in T_H}|K|\left(||b^{H}||_{0, \infty, K}||\nabla u_{H}||_{0, \infty, K}+||c^{H}||_{0, \infty, K} ||u_{H}||_{0, \infty, K}\right)||v_H-\Pi_H^*v_H||_{0, \infty, K}\\
\!\!\!&\leq &\!\!\!C\sum_{K\in T_H}\left(||b^{H}||_{0, \infty, K}||\nabla u_{H}||_{0, 2, K}+||c^{H}||_{0, \infty, K} ||u_{H}||_{0, 2, K}\right)||v_H-\Pi_H^*v_H||_{0, 2, K}\\
\!\!\!&\leq &\!\!\!C\left(||b^{H}||_{0, \infty, \Omega}+||c^{H}||_{0, \infty, \Omega} \right) ||u_{H}||_{1, \Omega}||v_H-\Pi_H^*v_H||_{0, \Omega}\\
\!\!\!&\leq &\!\!\!CH\left(||b^{H}||_{0, \infty, \Omega}+||c^{H}||_{0, \infty, \Omega} \right) ||u_{H}||_{1, \Omega}|v_H|_{1, \Omega}\\
\!\!\!&\leq &\!\!\!CH\left(||b^{0}||_{0, \infty, \Omega}+||c^{0}||_{0, \infty, \Omega}+e(HMM) \right) ||u_{H}||_{1, \Omega}|v_H|_{1, \Omega}\\
\!\!\!&\leq &\!\!\!C\left(H+e(HMM) \right) ||u_{H}||_{1, \Omega}|v_H|_{1, \Omega}\\
\!\!\!&\leq &\!\!\!C\left(H+\frac{\varepsilon}{\delta} \right) ||u_{H}||_{1, \Omega}|v_H|_{1, \Omega}.
\een
Therefore, by the estimate of $\varepsilon_1(u_H, v_H)$, $\varepsilon_2(u_H, v_H)$ and $\varepsilon_3(u_H, v_H)$, we obtain
\ben
|B(u_H, v_H)-A_{FVM}(u_H, v_H)|\leq C \left( H+\frac{\varepsilon}{\delta}\right)||u_H||_{1,~\Omega}||v_H||_{1,~\Omega}
\een\hfill\rule{1.5mm}{3mm}

In order to estimate the error between $u^0$ and $u^H$, we separate it into two parts
\ben
||u^0-u_H||_{1, \Omega}\!\!\!&\leq &\!\!\!||u^0-u^0_H||_{1, \Omega}+||u^0_H-u_H||_{1, \Omega},
\een
where $u^0_H$ is the numerical solution of (\ref{FEM}), which is FEM formula of the homogenized problem (\ref{homogenized problem}).
From the standard error estimate, we have
\be\label{FEM error}
||u^0-u^0_H||_{1, \Omega}\leq C H||u^0||_{1, \Omega}.
\ee
To estimate $||u^0_H-u_H||_{1, \Omega}$, the inf-sup condition for bilinear form $B$ in \cite{Wu} is useful.
\begin{lemma}\cite{Wu}
For $0<H<h_0$, we have
\ben
||u_H||_{1, \Omega}\leq C\sup_{0\neq v_H\in U_H}\frac{B(u_H, v_H)}{||v_H||_{1, \Omega}},~\forall u_H\in U_H.
\een
\end{lemma}
Then the second part can be derived as
\be\label{u0h_uh}
||u^0_H-u_H||_{1, \Omega}\!\!\!&\leq &\!\!\!C\sup_{0\neq v_H\in U_H}\frac{B(u^0_H-u_H, v_H)}{||v_H||_{1, \Omega}}\nonumber\\
\!\!\!&\leq &\!\!\!C\sup_{0\neq v_H\in U_H}\left(\frac{|B(u_H, v_H)-A_{FVM}(u_H, \Pi_H^* v_H)|}{||v_H||_{1, \Omega}}+\frac{|(f, v_H)_H-(f, \Pi_H^* v_H)_H|}{||v_H||_{1, \Omega}}\right)\nonumber\\
\!\!\!&:= &\!\!\!T_1+T_2.
\ee
By Lemma \ref{B_A}, we can easily obtain
\be\label{T1}
T_1\!\!\!&\leq &\!\!\!C\left(H+\frac{\epsilon}{\delta} \right) ||u_{H}||_{1, \Omega}\nonumber\\
\!\!\!&\leq &\!\!\!C\left(H+\frac{\epsilon}{\delta} \right) \left(||u^0||_{1, \Omega}+||u_{H}-u^0||_{1, \Omega}\right).
\ee
Denote by
\ben
\varepsilon(f)\!\!\!&= &\!\!\!(f, v_H)_H-(f, \Pi_H^* v_H)_H
=\sum_K\int_K f v_H d x - \sum_K |K| f(Q) \Pi_H^* v_H\\
\!\!\!&= &\!\!\!\sum_K\int_K f v_H d x - \sum_K \int_K f(Q) \Pi_H^* v_H dx=\sum_K\int_K f v_H d x - \sum_K \int_K f(Q)   v_H dx\\
\!\!\!&= &\!\!\!\sum_K\int_K (f(x)-f(Q))v_H dx \leq CH\sum_K||\nabla f||_{0, 2, K}||v_H||_{0, 2, K}\\
\!\!\!&\leq &\!\!\! C H |f|_{1, \Omega}||v_H||_{0, \Omega} \leq  C H |f|_{1, \Omega}||v_H||_{1, \Omega},
\een
then
\be\label{T2}
T_2\leq CH|f|_{1, \Omega}.
\ee

Therefor, from (\ref{FEM error}), (\ref{u0h_uh}), (\ref{T1}) and (\ref{T2}), we have
\ben
||u^0-u_H||_{1, \Omega}\leq C\left(H+\frac{\varepsilon}{\delta} \right).
\een\hfill\rule{1.5mm}{3mm}

\section*{Acknowledgments}
This work is supported in part by the Doctoral Starting up Foundation of Jinggangshan University under the Grant JZB11002.

\bibliographystyle{amsplain}

\end{document}